\documentclass{amsart}

\usepackage{amsmath,amsfonts,amsthm,amssymb,stmaryrd,paralist,tikz}
\usetikzlibrary{matrix,arrows,decorations.pathmorphing}

\usepackage{url,hyperref} % Required for links
\usepackage{graphicx} % Required for inserting images

\usepackage{color}
\definecolor{verydarkblue}{rgb}{0,0,0.4}
\usepackage{hyperref}
\hypersetup{
pdfauthor={Juliana Belding, Keaton Quinn},
pdftitle={Parameters in Calculus},
colorlinks=true,linkcolor=verydarkblue,
citecolor=verydarkblue,urlcolor=verydarkblue
}

\title{A Case for Teaching About Parameters in Calculus Courses}
\author{Juliana Belding}
\address{Department of Mathematics, Boston College, Chestnut Hill, MA, United States of America}
\email{beldingj@bc.edu}

\author{Keaton Quinn}
\address{Department of Mathematics and Statistics, Wellesley College, Wellesley, MA, United States of America}
\email{kq101@wellesley.edu}

\date{September 30, 2025}

\begin{document}

\begin{abstract}
Working with letters that represent unknown constants, i.e., parameters, has been historically challenging for students. This is an important skill for their success in many future quantitative settings, and yet it appears this topic is rarely included explicitly in math curricula. We argue that we should be explicitly teaching our students how to work with parameters, and that single variable calculus courses are a natural place to do so. We offer justification for this as well as examples and sample outlines for incorporating parameters into these classes. 
\end{abstract}

\maketitle

\section{Introduction} % Introduction and why do this in Calc 1 and in Calc  2

%In this paper, we describe our experiences including the concept of parameters in introductory calculus courses. We give the motivation for doing this, some key examples of parameters, and outlines of implementation in various courses. Lastly, we give some student and instructor feedback this. 

\subsection{Why Teach Parameters in a Math Course}

Consider a situation that has both variables and constants present, such as a model, function, equation, or expression. Sometimes the values of the constants are left unspecified, or we may imagine their values as being allowed to change somehow. We call these \emph{parameters}: any quantity that is considered constant for a particular situation, but can vary from one setting to the next \footnote{This definition can perhaps be mapped to the notion of the traces of a multivariate function. However there is merit in considering expressions with parameters as single variable functions with unspecified constants.}.
% most readily being that we don’t need to set up the formalism of a higher dimensional setting.}. 
For example, take the general formula for a quadratic function or parabola, 
\[
y = ax^2 + bx + c.
\]
The letters $a, b, c$ are parameters because they represent constants, while $x$ and $y$ are variables. Different choices of $a$, $b$, and $c$ (with $a \neq 0$) create different functions in the same ``family".

Consider the following models \footnote{These are, in order, the logistic model of population growth, a plasma concentration model from pharmokinetics, and a linear regression model from an econometrics textbook.} our students might encounter during a STEM education 
\[
P = \displaystyle {\frac {L}{1+Ae^{-kt}}}, \quad\quad
C = \frac{C_0e^{-k_et}}{1-e^{-k_e \tau}}, \quad\quad
y = \beta_0 + \beta_1x + u.
\]
Notice that these look different from the ``$y = f(x)$" mathematics they are more familiar with, and they may be intimidated by the multitude of letters. Working with parameters such as these is a skill that many students will need, whether in other disciplines or in higher-level math courses. In fact, the Mathematical Association of America 2004 report on the mathematical needs of partner disciplines explicitly mentions working with parameters as a desired skill in many subjects, including Biology, Health-related Life Sciences, Chemistry and Engineering \cite{ganter2004}. 

Furthermore, analyzing functions or situations that have parameters also helps students directly experience a key aspect of mathematics, the power and efficiency of generalizing. We can understand and describe the behavior of a whole class of functions versus doing it on a case-by-case basis. Take, for example, the logistic model mentioned above. We can compute the limit as $t \to \infty$ of the model with parameters to notice the numerator is always the population far in the future. Now we know how \emph{any} population that follows logistic growth will behave in the long run.

In our experience, students are often expected to use these ideas without being formally taught. We believe that we should be placing an explicit emphasis on teaching our students to work with parameters in our courses and that Calculus 1 and 2 are natural places to do this because 
\begin{itemize}
    \item there are interesting questions and results in Calculus that can be addressed using parameters (see examples in the next section); 
\item for many students this is their last math class and is thus a last chance to address this topic from a mathematical perspective;
\item for students who continue on to higher math and science courses, parameters will be used extensively and familiarity is often expected (e.g.,  differential equations, physics, or probability courses). 
\end{itemize}

\subsection{Challenges with the Concept of Parameters}

%(Q: Do we want to disambiguate parameters language? Give some note that this words means many things in different disciplines. See \url{https://en.wikipedia.org/wiki/Parameter}.

%At minimum, we should discuss what we take it to mean and that students have their own history (I've had students think of it as "limits" or boundaries of a problem). 

%Wikipedia: Statistical parameter a parameter is any quantity of a statistical population that summarizes or describes an aspect of the population, such as a mean or a standard deviation

The word parameter has many meanings \footnote{For example, in statisics, a  parameter is a quantity for the population that describes an aspect of the population, such as a mean or a standard deviation}. 
The definition given above follows a framework for types of variables proposed by Thompson, see \cite{thompson2017}. We note that this concept is itself challenging as it includes multiple ways to use parameters: as placeholders for unknown or unspecified constants, or as ``dynamic sliders” that can be changed to explore the effect of parameters on a situation \cite{drijvers2001}. 

Another challenge is that there are situations where there are choices for what should be considered variable versus parameter. For example, in the classic equation for the value of money with compound interest
$$M = Ae^{rt}$$
we have four quantities, $M, A, r, t$. We might want to explore the relationship of $M$ as a function of the interest rate $r$. In that case, $A$, the initial balance, and $t$, the time of investment are parameters, fixed but unspecified. In other cases we might be interested in $A$ as a function of $t$, if we have a fixed value of $M$ in mind. Determining what the parameters are thus depends on the context. 

This is something that students will need practice with in order to become fluent with these ideas. Bardini, Radford, and Sabena in \cite{bardini2005} found, while trying to probe students' understanding of this distinction, that their understanding of the idea of a variable was not as solid as expected. Park and Rizzolo \cite{park2022} noticed that even our teaching assistants might not have as good of a grasp on the idea of a parameter as we think. Yet Gray, Loud, and Sokolowski have determined that ``when students are attempting to learn the concepts $\ldots$ they must possess a robust and flexible view of all aspects of the changing quantities that form the foundation of their mathematical studies" \cite{gray2007}. It appears to us that there is a gap in mathematics curricula. These are difficult ideas that students apparently are expected to learn on their own, or simply from repeated exposure. Instead, we believe many could benefit from a more intentional treatment of parameters.

%Quote from Wikipedia We refer to the relations which supposedly describe a certain physical situation, as a model. Typically, a model consists of one or more equations. The quantities appearing in the equations we classify into variables and parameters. The distinction between these is not always clear cut, and it frequently depends on the context in which the variables appear. Usually a model is designed to explain the relationships that exist among quantities which can be measured independently in an experiment; these are the variables of the model. To formulate these relationships, however, one frequently introduces "constants" which stand for inherent properties of nature (or of the materials and equipment used in a given experiment). These are the parameters.[1]

\subsection{Why Parameters in the Intro Calculus Sequence}

As mentioned, working with parameters is a skill many students will need in other disciplines. 
For example, only a small subset of students will need to take the derivative or integral of a complicated function by hand in a later course or career. However, many of them may encounter models or formulas that contain parameters. It is thus helpful for students in math classes to develop strategies to be able to parse such models or formulas in other contexts.

In a first Calculus course (or even a pre-calculus course), there are many levels for students to work with parameters. Students can begin learning the concept in the context of functions and explore this with examples and dynamic technology. They can then learn to do computations with functions involving parameters, such as limits and derivatives. Lastly, they can learn to reason and interpret results that involve parameters, for example, how the limit at infinity of a logistic function depends on parameters.  
%[Goals(levels): 0: Understand concept of parameters, possibly with technology  1. Do mechanical computations with equations, etc involving parameters; 2. Careful interpretation of technical results]

% I've updated this paragraph. Will still give it some more though. Didn't know how explicit to go with examples since the next section is about examples. I LIKE THIS OVERALL. SHORT AND SWEET. 
A second calculus course is a natural time to place an such an emphasis, since they arise so frequently in the material. A general issue for students, in our opinion, is that the density of notation in Calculus 2 seems much higher than in a first course. Letters are much more common in the definitions and statements we present, such as for the definitions of definite integrals and infinite series, and students need to be prepared for the flexibility in their use. Being able to determine whether a letter represents a variable, a parameter, a constant, an index etc., is necessary for understanding even the basics of the content of the class. Indeed, they need to be able to parse these expressions to even begin a problem.

Moreover, much of the material deals with generalizing or classifying the objects in question. We shift away from discussing examples to discussing whole classes of examples. We build theory by preforming computations on families of objects by using parameters instead of explicit numbers. For example, the geometric series convergence test takes the form of checking whether the series you have is of the correct form, and then determining the value of a parameter related to the series (the common ratio). Similarly, the behavior of a solution to a differential equation, e.g., $y' = ky$, is in many cases governed by the values of the parameters in the equation, and we can make statements about the solutions without explicitly computing them.  Within Calculus 2, we expect students to be comfortable with the idea of a parameter, to have the mechanical skills to perform computations with them, and to have the reasoning skills to analyze or interpret their results.  Familiarity with parameters is almost expected by this point in the sequence, yet rarely do we dedicate time to teach them these skills.

We note that while these courses are a good place to talk about parameters, from our experience working with them is rarely included in a traditional syllabus for the course. We, perhaps unconsciously, set these learning goals without actually discussing them with students.  Many textbooks include problems in the exercises with parameters or results given in terms of parameters, but there is not an explicit section of the text dedicated to defining ``parameter" or concrete strategies to approach problems with parameters. In the next sections we give some examples of how we have tried to fill this gap.

\section{Examples} % Main demonstrative example from Calc 1 and 2 each
\label{sec:examples}

In this section, we share a few examples of problems involving parameters that we ask students to engage with in our calculus courses. One example is from Calculus 1 for Life and Social Sciences and the other is from Calculus 2 for Math and Science majors. A commonality among the problems here (and many such ``parameter problems") is that we start with analyzing a particular example and then transition to analyzing a generalization of the example. This is a classic mathematician move and allows students to experience the power of this first hand. 

A difference between the two examples is found in their level of abstraction. In the Calculus 1 example, the focus is on functions and modeling: parameters allow us to describe a more general function/model which makes the model more widely applicable. In the Calculus 2 example, parameters appear in a more abstract setting, an improper integral. They help generalize a result about such integrals that itself is then used as a tool in other contexts.

\subsection{Example from Calculus 1}

A classic example for a Calculus 1 class is the surge function. 
For example, here is a problem found in a standard textbook \cite[Problem 73, Section 4.1]{stewart2021}. 

{\it After the consumption of an alcoholic beverage, the concentration of alcohol in the bloodstream (blood alcohol concentration, or BAC) surges as the alcohol is absorbed, followed by a gradual decline as the alcohol is metabolized. The function
$$C(t) = 1.35te^{-2.802t}$$
models the average BAC, measured in mg/mL, of a group of eight male subjects $t$ hours after rapid consumption of 15 mL of ethanol (corresponding to one alcoholic drink). What is the maximum average BAC during the first 3 hours? When does it occur?}

Students can use derivatives to analyze the maximum value of this specific function, and find that the maximum average BAC is 0.18 mg/ml at $t = \frac{1}{2.802} \approx 0.35$ hours after consumption (rounded to nearest hundredth).  %\frac{1}{2.802} \approx 
They could also take the limit at infinity to see that the BAC will approach zero in the long run. 

After doing this, students can go further and look at a general surge function model $$C(t) = ate^{-kt}$$
for the concentration $C$ of a drug in the bloodstream at time $t$ \cite[Section 4.8]{hughes2021}.
They can ask what happens in general and how the parameters $a$ and $k$ affect the end behavior and maximum value. Students can explore this graphically with a dynamic graphing tool such as Desmos \footnote{See for example \url{https://www.desmos.com/calculator/abc9d0d0b3}.} to make some conjectures about how $a$ and $k$ affect the graph. Then, they can take the derivative of $C(t)$ and determine that the time of maximum concentration is $t = \frac{1}{k}$. Thus only $k$ affects when the peak occurs.

%They can also determine that the and that only $k$ affects the time of max concentration while both $a,k$ affect the max value. 
%they can realize that drug concentration approaches zero regardless of $a, k$. 

As an extension, they can answer interesting questions such as ``what values of $a$ and $k$ will guarantee that the drug concentration does not exceed a certain threshold for safety?". Students can again explore this with a graphing tool then analytically solve for a relationship of $a,k$ that guarantees this.

Note that to get students to this point, they need to first understand what a parameter is conceptually and then develop practice taking derivatives and finding critical points of functions with parameters.

\subsection{Example from Calculus 2}

% We should discuss the style of this, if it's formal enough, etc. Also if I go too detailed on solving the problem. 

%Juliana - I like the tone and I think I could make my section more conversational. We could even explicitly say that we are writing it that way

Perhaps the prototypical example from a second course in calculus is the $p$-integral. When discussing improper integration, popular examples are $\int_1^\infty \frac{1}{x} dx$ and $\int_1^\infty \frac{1}{x^2} dx$ since they are so similar yet have opposite behavior. Depending on the course, one might take this a step further. We pose the following problem to students. 
 
\emph{
For what value(s) of $p$ does the following integral converge? 
\[
\int_1^\infty \frac{1}{x^p} dx
\]
}

\noindent This family is the perfect illustration since investigating its behavior requires both technical and analytical skills to determine. 

The lecture usually proceeds as follows:
\emph{\begin{itemize}
\item We can investigate convergence directly by integrating. How would we integrate this? 
\item Use the power rule.
\item Can we always do that? Are there any values of $p$ where that wouldn’t work? 
\item Oh, first we need to know if $p = 1$ or not. 
\item Ok, if $p=1$ what does the integral do? It diverges. If $p\neq 1$ then we can use the power rule to get
\[
\int_1^\infty \frac{1}{x^p} dx = \lim_{t \to \infty} \frac{1}{1-p}\left[ t^{1-p} - 1\right],
\]
and when does this converge? 
\item When the limit exists. Which means we need the exponent on $t$ to be negative $1-p < 0$. So it must be that $p > 1$.
\end{itemize}}
\noindent
In order to do all of this, students must be able to perform the integration with the parameter, and then also reason about what the values of $p$ imply about convergence. 

To illustrate why this ability is so central to Calculus 2, consider another example: 

\emph{Does $\sum \frac{1}{\sqrt{n^4 - 1}}$ converge? 
\begin{itemize}
\item Let’s make a comparison to $\sum \frac{1}{n^2}$. What does this do?
\item Well, let's compare this to the integral $\int_1^\infty \frac{1}{x^2} dx$, and this converges because $p = 2 > 1$. 
\end{itemize}}
Here students begin to see the power of abstraction. We can now skip direct computation as we have built some theory. If we understand the behavior of the $p$-integrals, then we can answer other questions in broader settings. And we do understand the behavior simply from the value of the parameter $p$.

\section{Implementation} % Our specific implementations in our classes, as a guide for others

In this section, we describe two cases of implementing parameters in Calculus courses. 
For each, we describe the course context and audience, the learning goals, assessment, and an outline of course content related to parameters.
For more details, please contact the authors.

\subsection{Calculus 1: A Sample Outline}

%Description of course/class setting 
%Organize as overall learning goals 
%unit by unit with some examples 
%Assessment - formative and summative 
%Student feedback
%46 of about 400 self-report as new to calc

\subsection*{Course Audience and Context} The course is a Calculus 1 course for life and social science majors taught between 2022-2024 at Boston College. The content is limits and derivatives through optimization, with a variety of functions including exponential, logarithmic, and trigonometric functions. The course is a multi-section coordinated course of about 400 students each fall. Instructors include experienced professors as well as post-docs and graduate student instructors. Students have a MWF class led by their instructor and a T or TH discussion section with a graduate teaching assistant. We have included this material for about three years with positive results overall. 

We note  another reason for including this material in the course, in addition to those stated in Section 1.1. Close to 90\% of students in this course have had Calculus before, and we included parameters as a topic that would both challenge students and be relevant in their other courses (versus increasing challenge through more advanced algebra or trigonometry required for problems, or including more theoretical topics like Squeeze Theorem).

\subsection*{Major Learning Goals for Parameters} We want students to be able to do the following by
the end of the course.
\begin{itemize}
    \item Given an equation or function or model, students can identify parameters and variables from the surrounding text and notation, and identify underlying structure in simple cases like linear or exponential. 
    \item Students are able to use limits and derivatives to analyze features of a model with parameters, such as end behavior and extrema. 
\end{itemize}

\subsection*{Assessment} Problems with parameters are included in all types of assessment.  There is usually about one such problem weekly in written homework and every other week in discussion section problems. They also appear occasionally in online homework. 
 
There are three midterms and a final. A problem with parameters and limits appears on the second midterm.  Problems involving computation of derivatives for a simple function with parameters like $y  = \cos(at)$ may appear on the third midterm. A question with parameters in the context of rates of change and/or optimization appears on the final exam. 

\subsection*{Course Outline and Teaching Notes} Here is an outline of how the topic of parameters  appears in the course, which is 15 weeks. 
Below are the learning goals for each section and some brief descriptions. 
\begin{itemize}
    \item Introduction to the concept of parameter and dynamic graphing (Week 2) 
    \item Parameters in piecewise functions (Week 5) 
    \item Limits  of functions with parameters (Week 6)
    \item Derivatives of functions with parameters (Week 10) 
    \item Max and min of functions with parameters (Week 15)
    \end{itemize}

{\bf Introduction to the concept of parameter and dynamic graphing.} In the second week of the semester, we introduce the notion of a parameter as an unspecified constant. 
We discuss how parameters allow us to describe a whole family of functions or equations that have similar features, and we work with familiar examples like linear, quadratic and  exponential functions. We discuss how parameters are often directly related to key features, such as $m$, the slope, of a linear function $y = f(x) = mx+b$.

The goals for students are two-fold. First, students can graph a function with parameters in Desmos, namely how to add ``sliders" for the parameters and set the range for parameters. They are able to qualitatively describe what happens to the function and its key features when we change parameters. Secondly, given an equation or function that includes parameters and variables, students can identify which is which from the surrounding text and notation. They can determine the function's structure (eg: linear, exponential, inverse proportional, etc.) by picking sample parameter values of appropriate sign.
   
           For example, consider the problem below. 

\textit{
Consider someone who is a regular coffee drinker. Their demand for cups of coffee, $Q$, depends on the price of a cup of coffee, $P$. 
One possible model is
$$Q = \frac{ \alpha M}{P}$$
where $P$ is in dollars, $M$ is the weekly income (in dollars) of the consumer and $\alpha$ is the fraction of their income they spend on coffee. 
}

Students can pick a sample value $\alpha = 0.05$ and $M = 1200$ and see that $Q(P)$ is an inversely proportional function.

%(eg: we can often determine the meaning/impact of the parameter). 

{\bf Parameters in piecewise functions.} When we learn about continuity and differentiablity, students work with parameters in piecewise functions. The goal is for students to determine the parameter values in piecewise functions that make a function continuous and/or differentiable at a point. In some situations, we supplement the analytic solving process with a Desmos exploration so they can better visualize what is going on with the pieces of the function.

%tudents solve classic problems where they determine the values of parameters in piecewise functions that make a function continuous and/or differentiable at a point.

{\bf Limits  of functions with parameters.} At the end of the unit on computing limits, students compute limits for functions with parameters. The focus is primarily on limits at infinity. They begin to understand that the answer may or may not involve the parameter(s) and this gives us information about which parameters affect the end behavior of a function.

%The goals for students are as follows: 
%\begin{itemize}
%\item Students can compute derivatives of functions with parameters. 
%\item Students can graphing tools to qualitatively describe the effect of parameters on extrema. 
%\item Students can use derivative tools to determine exact location of extrema in terms of parameters. 
%\end{itemize}
{\bf Derivatives and Extrema of functions with parameters.} At the end of the unit on computing derivatives, students compute derivatives of functions with parameters. They can then identify which parameters affect the derivative of a function at a specific point, like the origin. They can confirm this visually using Desmos and changing the parameter and looking at the slope of the graph. 

Lastly, we discuss parameters and optimization. Students consider how parameters affect the maxima or minima of a function. They begin with graphing the function with parameters using Desmos, qualitatively describing the effect of the parameter on the $x$ and $y$-coordinates of an extremum. Then they use derivative tools to determine the exact locations and answer questions such as the one about surge functions in the Examples section above.   
In this section, we give students the function to be optimized, as there is enough cognitive load without also having them create a function that models the situation (as in the drug concentration example).
 
We also highlight that parameters come up naturally in the ``classic textbook" optimization problems.
For example, after solving an optimization problem of minimizing surface area of a cylindrical can for a fixed volume of 355 cubic cm, we can ask ``what if the manufacturer changes their mind and wants to know the best choice for a volume of 400 cubic cm?". Students realize  they can solve the optimization problem with a generic volume $V$ as parameter. Many students seem to find it natural to think about parameters here and are motivated by the idea that we can be efficient and find the optimal radius and height for any possible $V$ value.

\subsection{Calculus 2: A Sample Outline}

\subsection*{Course Audience and Context}
Our Calculus 2 classes were taught in 2023 at Boston College for an audience of Math and Physical Sciences majors, and in 2024-25 at Wellesley College for a general audience, including math majors. These were semester courses consisting of one or two sections with 30 students each, and taught by a single instructor. We focused on integration and infinite series, starting with Riemann sums and ending with Taylor series.  

\subsection*{Major Learning Goals for Parameters} 
We wanted students to be able to do the following by the end of the course. 
\begin{itemize}

\item Students should become comfortable working with expressions and equations with unspecified constants. 

\item Students should be able to recognize an expression with parameters as a whole family of expressions, and understand how the example moves within the family as the values of the parameters change. 

\item Students should be able to answer questions by referring to the structure of the setting and the values of any parameters, as in the $p$-integral example, without performing any further computations. 

\end{itemize}

\subsection*{Assessment} 
Problems with parameters were included in all types of assessment. Student learning was assessed using a free and open-source online homework system, using written homework assignments, and on 3 midterm exams and a final exam. Online homework was assigned every lecture period and one problem per set was usually a standard computation with parameters present. 
On written homework one problem per set asked them to work more deeply with parameters by, for example, finding how their values affected behavior of the functions present. 

\subsection*{Content Outline and Teaching Notes}

We start with introducing the notion of a parameter in the first week of the semester. We do examples of determining which letters in the equations are the independent and dependent variables, and which are parameters and constants. Then we discuss more examples of how we can interpret the same equation as different types of functions depending on what we consider the independent variable. This helps prime them for area functions where we change the endpoint of an interval of integration from a parameter to a variable. Other examples serve as a review of some Calculus 1 ideas, since students probably need refreshing after the break between semesters. These examples include determining which parameters in equations affect the derivative of a function, or the limit of the function, or the location/value of any extrema.

{\bf Riemann Sums and Definite Integrals.} 
We move into definite integrals and the FTC. This setting is ripe with parameters. They show up as the endpoints of intervals, as the number of rectangles in a Riemann sum, or as constants of integration. Desmos is useful here as students can adjust values to see how the geometric picture of integration is affected. Starting with the first midterm, each summative assessment has a question dealing with parameters, typically a computation problem asking students to determine an integral or limit that contains letters. Written homework is another opportunity to for students to practice these ideas. 

{\bf Integration Techniques and Applications.}
The class moves on to integration techniques and applications. We begin to sprinkle parameters into the lecture questions and homework sets. “How would we integrate this function if this number was a letter?”. Or, Integrate $\int \frac{1}{ax + b} dx$ for $a \neq 0$. Or, find the antiderivative of  $\frac{1}{x^2 + a^2}$. Students see that having antiderivatives of families is useful when, e.g., integrating rational functions using partial fractions. Here we can decompose a complicated integral into sums of a few known types, and we can save time finishing the problem by appealing to their answers to these homework questions. Written homework is also an opportunity to include examples from other disciplines. For example, one can solve a reaction rate equation with parameters from chemistry by turning the separable differential equation into an integral equation. After discussing the $p$-integrals, we begin having students investigate the effect parameters have on convergence. Questions like ``for what value(s) of $k$ does the integral $\int_1^\infty \frac{x}{x^2+1} - \frac{k}{x+1} dx$ converge?” become more common.

{\bf Infinite Sequences and Series.}
Familiarity with parameters defining families becomes almost necessary when reaching infinite sequences and series. These are more challenging to view geometrically using Desmos, for example. It is harder to visualize the affect their values have, so maturity with them is helpful. But here is where the idea hits the hardest, as most convergence tests deal with classifying series into nice buckets and then comparing wild series to nice ones. These ideas occur naturally in this material, and their incorporation  in the class starts to feel less contrived. 

{\bf Power and Taylor Series.}
Power series and Taylor series round out the semester. The center of a series and its radius of convergence have focus. More classification problems come up in the form of determining if a function is analytic. We end with applications of Taylor series from other disciplines. For example, we have students find the second order Taylor approximation to the function $f(x) = \frac{1}{\sqrt{1-x^2}}$ at $x = 0$ and then introduce total relativistic energy as a modification of this function with parameters, producing
\[
E(v) \approx mc^2 + \frac{1}{2}mv^2,
\]
for small velocities $v$. This usually captures their interest as they recognize both the iconic $mc^2$ from Einstein and $1/2mv^2$ as kinetic energy from their physics classes.

\section{Conclusion}
We believe that we can and should explicitly teach our students how to work with parameters. This can  strengthen their ability to make sense of mathematical models and formulas they may encounter in future courses. Furthermore, much of the content of their classes is on describing or classifying general behavior, and parameters are a language to do so. Yet it appears to us that students are often expected to pick things up as they go along, with varied success. Our solution to this is to intentionally incorporate the topic into single-variable Calculus courses.

Going forward, we are interested in examining how effective our curriculum is in helping students learn to work with parameters and measuring their understanding, both directly with exams and indirectly via students self-report. We have some small-scale positive results about this.  Students in the first author's Calculus I section were asked about their agreement with the following statement:
 {\it
``I feel more confident working with a function with parameters compared to when I started the course."
}
The students \footnote{Of the 101 students in the section,  78 responded, for a 77\% response rate.}  responded overwhelmingly positively: 87\% of students agreed or strongly agreed, with 12\% neutral and 1\% disagreeing.

For future directions, we suggest considering where and how parameters might be explicitly addressed in the context of other courses at the high school or college level, such as in college algebra or precalculus. More generally, we imagine there are other topics in mathematics that instructors expect students to become comfortable with that are not explicitly addressed. We encourage instructors to reflect on such topics in their own courses.

\section*{Acknowledgments}
We would like to thank Oscar Fernandez, J\=anis Lazovskis, and Erika Ward for helpful comments on drafts of this paper.

%\bibliography{References}

\begin{thebibliography}{HHGLF21}

\bibitem[BRS05]{bardini2005}
Caroline Bardini, Luis Radford, and Cristina Sabena.
\newblock Struggling with variables, parameters, and indeterminate objects or
  how to go insane in mathematics.
\newblock In {\em Proceedings of the 29th Conference of the International Group
  for the Psychology of Mathematics Education}, volume~2, pages 129--136, 2005.

\bibitem[Dri01]{drijvers2001}
Paul Drijvers.
\newblock The concept of parameter in a computer algebra environment.
\newblock In {\em Proceedings of the 25th Annual Meeting of the International
  Group for the Psychology of Mathematics Education}, volume~2, pages 385--392,
  Utrecht: The Netherlands, 2001.

\bibitem[GB04]{ganter2004}
Susan~L Ganter and William Barker.
\newblock {\em The curriculum foundations project: Voices of the partner
  disciplines}.
\newblock Mathematical Association of America, Washington, DC, 2004.

\bibitem[GLS07]{gray2007}
Susan~S Gray, Barbara~J Loud, and Carole Sokolowski.
\newblock Calculus students' difficulties in using variables as changing
  quantities.
\newblock In {\em Electronic Proceedings for the Tenth Special Interest Group
  of the Mathematical Association of America on Research in Undergraduate
  Mathematics Education: Conference on Research in Undergraduate Mathematics
  Education}, pages 1--15, San Diego, CA, 2007. Mathematical Association of
  America.

\bibitem[HHGLF21]{hughes2021}
Deborah Hughes-Hallett, Andrew~M Gleason, Patti~Frazer Lock, and Daniel~E
  Flath.
\newblock {\em Applied calculus}.
\newblock John Wiley \& Sons, 2021.

\bibitem[PR22]{park2022}
Jungeun Park and Douglas Rizzolo.
\newblock Use of variables in calculus class: focusing on teaching assistants'
  discussion of variables.
\newblock {\em International Journal of Mathematical Education in Science and
  Technology}, 53(1):165--189, 2022.

\bibitem[SCW21]{stewart2021}
James Stewart, Dan Clegg, and Saleem Watson.
\newblock {\em Single variable calculus: early transcendentals}.
\newblock Cengage, 9th edition, 2021.

\bibitem[TC17]{thompson2017}
Patrick~W Thompson and Marilyn~P Carlson.
\newblock Variation, covariation, and functions: Foundational ways of thinking
  mathematically.
\newblock {\em Compendium for research in mathematics education}, pages
  421--456, 2017.

\end{thebibliography}
%\bibliographystyle{alpha}

\end{document}